\documentclass[journal]{IEEEtran}
\ifCLASSINFOpdf
   \usepackage[pdftex]{graphicx}
\else
\fi
\usepackage[cmex10]{amsmath}
\usepackage[tight,footnotesize]{subfigure}

\usepackage[caption=false,font=footnotesize]{subfig}
\providecommand{\abs}[1]{\left|#1\right|}
\begin{document}
%
\title{Normalizations of the Proposal Density\\in Markov Chain Monte Carlo Algorithms}
%
%
%

\author{{Antoine~E.~Zambelli}
\thanks{This is a preprint. Official peer-reviewed version set to appear as \emph{Invited Paper} in \emph{Proceedings of MMSSE 2015}.}}

%
%

\markboth{Journal of \LaTeX\ Class Files,~Vol.~6, No.~1, January~2007}%
{Shell \MakeLowercase{\textit{et al.}}: Bare Demo of IEEEtran.cls for Journals}
%



\maketitle
\thispagestyle{empty}

\begin{abstract}
We explore the effects of normalizing the proposal density in Markov Chain Monte Carlo algorithms in the context of reconstructing the conductivity term $K$ in the $2$-dimensional heat equation, given temperatures at the boundary points, $d$. We approach this nonlinear inverse problem by implementing a Metropolis-Hastings Markov Chain Monte Carlo algorithm. Markov Chains produce a probability distribution of possible solutions conditional on the observed data. We generate a candidate solution $K'$ and solve the forward problem, obtaining $d'$. At step $n$, with some  probability $\alpha$, we set $K_{n+1}=K'$. We identify certain issues with this construction, stemming from large and fluctuating values of our data terms. Using this framework, we develop normalization terms $z_0,z$ and parameters $\lambda$ that preserve the inherently sparse information at our disposal. We examine the results of this variant of the MCMC algorithm on the reconstructions of several $2$-dimensional conductivity functions.
\end{abstract}


%
\IEEEpeerreviewmaketitle

\section{Introduction}
%
%
%
%
\IEEEPARstart{T}{he} idea of an inverse problem is to reconstruct, or retrieve, information from a set of measurements. In many problems, the quantities we measure are different from the ones we wish to study; and this set of \emph{d} measurements may depend on several elements. Our goal is thus to reconstruct, from the data, that which we wanted to study. In essence, given an effect, what is the cause? For example: If you have measurements of the temperature on a surface, you may want to find the coefficient in the heat equation.\\

The nonlinearity and ill-posedness of this problem lends itself well to Markov Chain Monte Carlo algorithms. We detail this algorithm in later sections, but we note now that there has been much work done on Metropolis-Hastings MCMC algorithms. However, much of it has been trying to determine optimal proposal densities (\cite{luengo},\cite{rosenthal}). Luengo and Martino (\cite{luengo}) treat this idea by defining an adaptive proposal density under the framework of Gaussian mixtures. Our work, however, is focused on improving the reconstruction given a proposal density.\\

We take no views on the optimality of the structure of the proposal density in our case, which we take from~\cite{fox}. We simply observe possible improvements to this density by normalizing it's terms through context-independent formulations. Eventually, we would like to implement the GM-MH algorithm of~\cite{luengo} on our proposal density, and provide a rigorous definition of our construction in an analogous manner to their work.\\

The paper is structured as follows. We first present the framework of our problem in the subsection below. Section~\ref{sec:mhmcmc} presents the MHMCMC algorithm and proposal densities along with non-normalized results. The error analysis of those results (in Section~\ref{sec:error}) motivates this work while Sections~\ref{sec:prelimnorm} to~\ref{sec:localnorm} present the new constructions and associated results.

\subsection{Heat Diffusion}

In this problem, we attempt to reconstruct the conductivity $K$ in a steady state heat equation of the cooling fin on a CPU. The heat is dissipated both by conduction along the fin and by convection with the air, which gives rise to our equation:
\begin{equation}\label{eq:heatpde}
u_{xx}+u_{yy}=\frac{2H}{K\delta}u
\end{equation}
with $H$ for convection, $K$ for conductivity, $\delta$ for thickness and $u$ for temperature. The CPU is connected to the cooling fin along the bottom half of the left edge of the fin. We use standard Robin Boundary Conditions with
\begin{equation}\label{eq:robinbc}
Ku_{normal}=Hu
\end{equation}
Our data in this problem is the set of boundary points of the solution to (\ref{eq:heatpde}), which we compute using a standard Crank-Nicolson scheme for an $n \times m$ mesh (here $20 \times 20$). We denote the correct value of $K$ by $K_{\textrm{correct}}$ and the data by $d$. In order to reconstruct $K_{\textrm{correct}}$, we will take a guess $K'$, solve the forward problem using $K'$, obtaining $d'$, and compare those boundary points to $d$ by implementing the Metropolis-Hastings Markov Chain Monte Carlo algorithm (or MHMCMC).

\section{Metropolis-Hastings MCMC}
\label{sec:mhmcmc}

Markov Chains produce a probability distribution of possible solutions (in this case conductivities) that are most likely given the observed data (the probability of reaching the next step in the chain is entirely determined by the current step). The algorithm is as follows (see~\cite{fox}). Given $K_n$, $K_{n+1}$ can be found using the following:
\begin{enumerate}
	\item Generate a candidate state $K'$ from $K_n$ with some distribution $g(K'|K_n)$. We can pick any $g(K'|K_n)$ so long as it satisfies
	\begin{enumerate}
		\item $g(K'|K_n)=0 \Rightarrow g(K_n|K')=0$
		\item $g(K'|K_n)$ is the transition matrix of the Markov Chain on the state space containing $K_n,K'$.
	\end{enumerate}
	\item With probability
	\begin{equation}\label{eq:alpha}
	\alpha(K'|K_n)\equiv min\left\{1,\frac{Pr(K'|d)g(K_n|K')}{Pr(K_n|d)g(K'|K_n)}\right\}
	\end{equation}
	set $K_{n+1}=K'$, otherwise set $K_{n+1}=K_n$ (ie. accept or reject $K'$). Proceed to the next iteration.
\end{enumerate}
More formally, if $\alpha > u\sim U[0,1]$, then $K_{n+1}=K'$. Using the probability distributions of our example, (\ref{eq:alpha}) becomes
\begin{multline}\label{eq:alpha2}
\alpha (K'|K_n)\equiv\\
\min\left\{ 1,e^{\frac{-1}{2 \sigma^2}\sum_{i,j=1}^{n,m} \left[ \left( d_{ij}-d_{ij}' \right)^2 - \left( d_{ij}-d_{n_{ij}} \right)^2  \right] }\right\}
\end{multline}
where $d'$ and $d_n$ denote the set of boundary temperatures from $K'$ and $K_n$ respectively, and $\sigma=0.1$. To simplify (\ref{eq:alpha2}), collect the constants and separate the terms relating to $K'$ and $K_n$:
\begin{align*}
\frac{-1}{2 \sigma^2}\sum_{i,j=1}^{n,m}{\left[ \left( d_{ij}-d_{ij}' \right)^2 - \left( d_{ij}-d_{n_{ij}} \right)^2  \right]}&=\frac{-1}{2}\left[ f' - f_n  \right]\\
&= -(D_1)
\end{align*}
Now, (\ref{eq:alpha2}) reads
\begin{equation} \label{eq:alpha3}
\alpha (K'|K_n) \equiv min\left\{ 1, e^{ -D_1  }\right\}
\end{equation}
Note that we are taking this formulation as given, and that the literature mentioned above (most notably Gaussian Mixture based algorithms) would be going from~\eqref{eq:alpha} to~\eqref{eq:alpha2} perhaps differently.

\subsection{Generating $K'$}

To generate our candidate states, we will perturb $K_n$ by a uniform random number $\omega\in[-0.005,0.005]$. In the simplest case, where we are dealing with a constant $K_{\textrm{correct}}$, then we could proceed by changing every point in the mesh by $\omega$, and the algorithm converges rapidly.\\

Looking at non-constant conductivities forces us to change our approach. If we simply choose to change one randomly chosen point at a time, then we have a systemic issue with the boundary points, which exhibit odd behavior and hardly change value. To sidestep this, we will change a randomly chosen grid ($2\times 2$) of the mesh at once. Thereby pairing up the troublesome boundary points with the well-behaved inner points.

\subsection{Priors}

While a gridwise change enables us to tackle non-constant conductivities, two issues remain. The first is that our reconstructions are still marred with ``spikes" of instability. The second, more profound, is that the ill-posedness of the problem means there are in fact infinitely many solutions, and we must isolate the correct one. This brings us to the notion of priors. These can be thought of as weak constraints imposed on our reconstructions. However, we do not wish to rule out any possibilities, keeping our bias to a minimum. So we define
\begin{multline}
T' =\sum_{j=1}^{n}\sum_{i=2}^{m} \left( K'(i,j)-K'(i-1,j) \right)^2\\ + \sum_{i=1}^{m}\sum_{j=2}^{n} \left( K'(i,j)-K'(i,j-1) \right)^2
\end{multline}
\begin{multline}
T_n =\sum_{j=1}^{n}\sum_{i=2}^{m} \left( K_n(i,j)-K_n(i-1,j) \right)^2\\ + \sum_{i=1}^{m}\sum_{j=2}^{n} \left( K_n(i,j)-K_n(i,j-1) \right)^2
\end{multline}
let $D_2=T'-T_n$, and modifying~\eqref{eq:alpha3}, we obtain
\begin{equation}\label{eq:alphac}
\alpha_c (K'|K_n) \equiv min\left\{ 1, e^{-\lambda_1 D_1 -\lambda_2 D_2} \right\} 
\end{equation}
By comparing the smoothness of $K'$ not in an absolute sense, but relative to the last accepted guess, we hope to keep as many solutions as possible open to us, while ensuring a fairly smooth result. We introduce one additional prior, this time imposing a condition on the gradient of our conductivity. The author explores the notion of priors more fully in~\cite{zambelli}, but much as we take the proposal density as given, the aim of this paper is not to examine priors per se. So we look at the mixed partial derivative of our candidate state and compare it to that of the last accepted guess
\begin{multline}
M'=\sum_{j=1}^{n}\sum_{i=2}^{m} \left( K_{xy}'(i,j)-K_{xy}'(i-1,j) \right)^2\\
+\sum_{i=1}^{m}\sum_{j=2}^{n} \left( K_{xy}'(i,j)-K_{xy}'(i,j-1) \right)^2
\end{multline}
\begin{multline}
M_n=\sum_{j=1}^{n}\sum_{i=2}^{m} \left( K_{n_{xy}}(i,j)-K_{n_{xy}}(i-1,j) \right)^2\\
+ \sum_{i=1}^{m}\sum_{j=2}^{n} \left( K_{n_{xy}}(i,j)-K_{n_{xy}}(i,j-1) \right)^2
\end{multline}
where $K_{xy}'$ and $K_{n_{xy}}$ are computed using central and forward/backward finite difference schemes. We let $D_3=M'-M_n$ and modify~\eqref{eq:alpha3} to get
\begin{equation}\label{eq:alphas}
\alpha_s(K'\mid K_n)\equiv\min\left\{1, e^{-\lambda_1D_1-\lambda_3D_3} \right\}
\end{equation}
We now take the acceptance step of our algorithm as
\begin{equation}\label{eq:alphacs}
\alpha=\max\left\{ \alpha_c,\alpha_s \right\}
\end{equation}
So the algorithm seeks to satisfy at least one of our conditions, though not necessarily both. We present some preliminary results in Figure~\ref{fig:tpfirst} and Figure~\ref{fig:gaussfirst} below. Note that we are clearly on the right path, with the algorithm approaching it's mark, but not to a satisfying degree.
\begin{figure}[h!]
	\centering
	\subfigure[Target.]{
		\includegraphics[height=3.1cm]{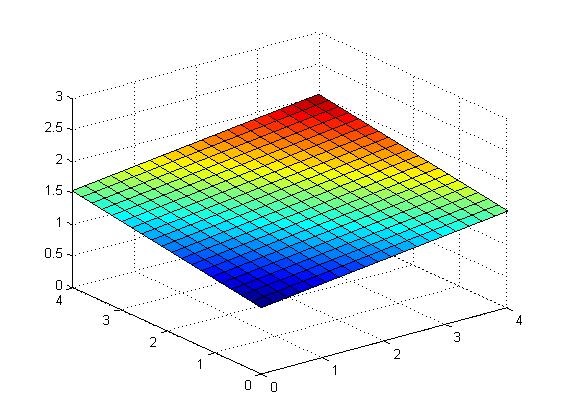}
	}
	\subfigure[Reconstruction with\newline $\lambda_1=1,\ \lambda_2=100,\ \lambda_3=15$.]{
		\includegraphics[height=3.1cm]{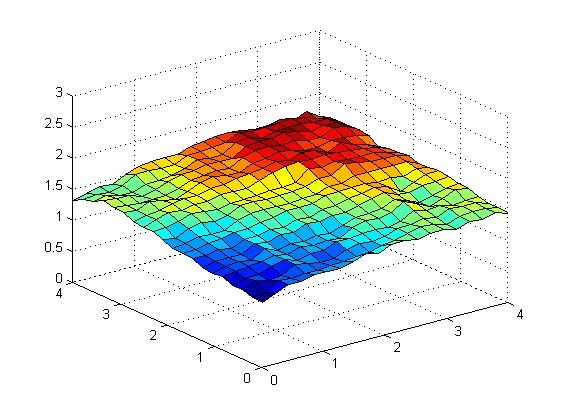}}
	\caption{Reconstruction of a tilted plane with priors, $10$ million iterations.}
	\label{fig:tpfirst}
\end{figure}
\begin{figure}[h!]
	\centering
	\subfigure[Target.]{
		\includegraphics[height=3.1cm]{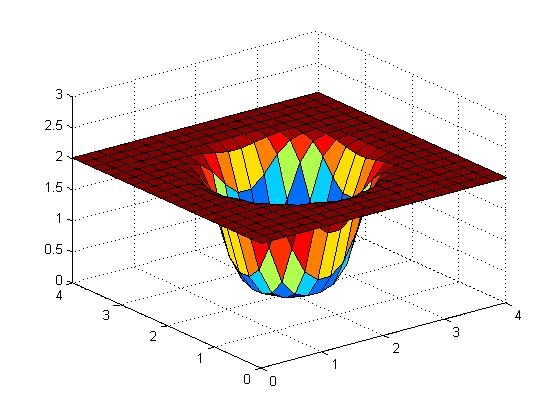}}
	\subfigure[Reconstruction with\newline$\lambda_1=1,\ \lambda_2=10,\ \lambda_3=15$ .]{
		\includegraphics[height=3.1cm]{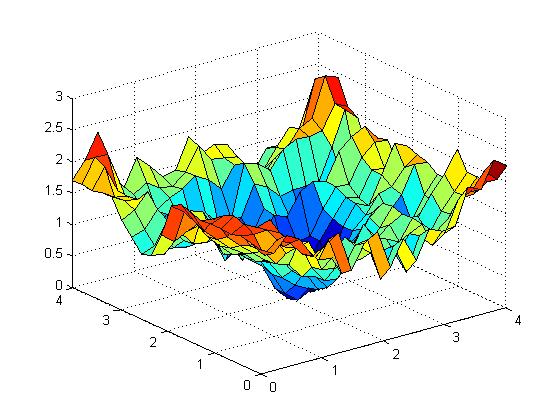}}
	\caption{Reconstruction of a Gaussian well with priors, $10$ million iterations.}
	\label{fig:gaussfirst}
\end{figure}

\section{Error Analysis}
\label{sec:error}
Our work so far has looked at qualitative improvements to our reconstructions, now we seek to quantify those improvements and the performance of the algorithm in general. Several metrics can be used for this purpose, but we will focus our writeup on the following: the difference between the data and the output using our guess ($\delta$), given by
\begin{equation*}
\delta=\left( \delta_1 \ \cdots \ \delta_n \right) \quad \textrm{, with } \delta_i=\sum{(d-d_i')^2} 
\end{equation*}
the sum of differences squared between $K_{correct}$ and $K_n$ ($\beta$),
\begin{equation*}
\beta=\left( \sum{(K_{correct}-K_1)^2}\  \cdots \ \sum{(K_{correct}-K_n)^2}\right)
\end{equation*}
and most importantly, the rate of acceptance of guesses ($\Gamma$), where
\[ \Gamma_0=0\quad \textrm{and} \quad
\Gamma_i =
\begin{cases}
\Gamma_{i-1}+1 & \text{if guess is accepted.}\\
\Gamma_{i-1}  & \text{if guess not accepted.}
\end{cases}
\]
for each subsequent iteration.\\

The form of $\Gamma$ is a step function, where accepting every guess would resemble a straight line of slope $1$, and accepting none of the guesses results in a slope of $0$. The shape of this function should tell us something about when the algorithm is performing best.

\subsection{$\delta$, $\beta$, $\Gamma$ Results}

\begin{figure}[h!]
	\centering
	\subfigure[$\beta$.]{
		\includegraphics[height=3cm]{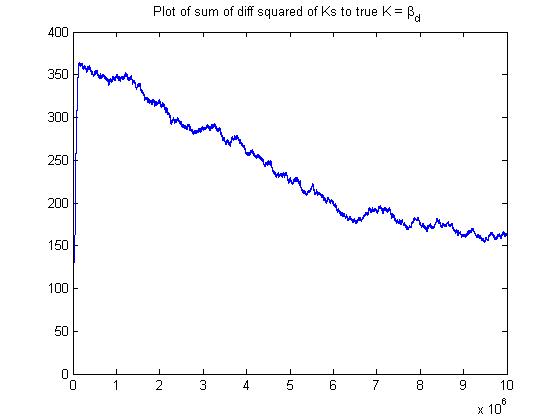}}
	\subfigure[$\delta$.]{
		\includegraphics[height=3cm]{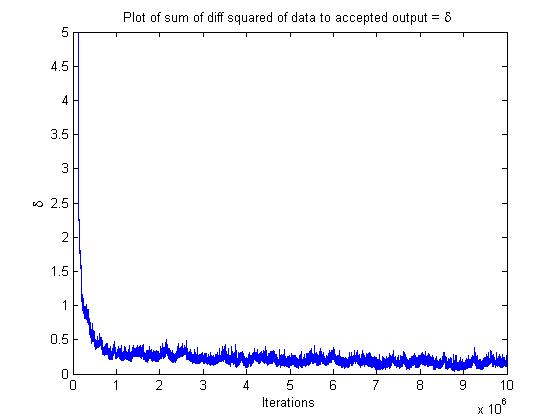}}
	\subfigure[$\Gamma$.]{
		\includegraphics[height=3cm]{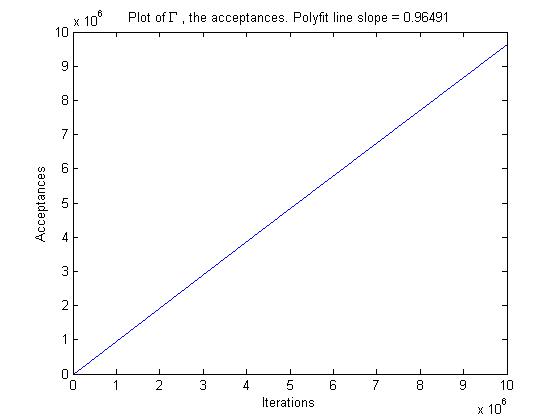}}
	\caption{Plot of the error metrics without normalizations.}
	\label{fig:errorimg}
\end{figure}

The results of tests involving these parameters reveals some interesting information (see Figure~\ref{fig:errorimg}). $\beta$ decreases, as expected, at a decreasing rate over time, slowing down around $6-7$ million iterations, which seems in line with the qualitative results.\\

On the other hand, $\delta$ decreases much more rapidly. The difference between the data and simulated temperatures becomes very small starting at as early as $250000$ iterations. In a sense, this fits with the problem of ill-posedness, the data is only useful to a certain degree, and it will take much more to converge to a solution (and we have been converging beyond $250k$ iterations).
The most important result, however, comes from $\Gamma$. If we fit a line to our step function, we get slopes of $0.95$ or more. This means we are accepting nearly every guess. While this could be troubling on its own, the fact that we are accepting at a constant rate as well is indicative of a deeper problem in our method.\\

Given that $\Gamma$ is dependent solely on the likelihood of accepting a guess, we take a look at $\alpha$ directly. What we find is that $\alpha$ is evaluating at $1$ almost every iteration. The quantities we are looking at within it (comparing data and smoothness) are simply too large. We need to normalize our distribution.

\section{Preliminary Structure}
\label{sec:prelimnorm}

In the following sections, we examine the impact of normalizations on our data terms, and explore the motivations behind the various constructions. More rigorous data is provided concerning the final form, while the earlier results focus on the concepts that guided their evolution.\\

One structural change which we will implement is to take equation~\eqref{eq:alphacs}, and change it to be more restrictive. Previously, it was looking for solutions which satisfied at least one of the prior conditions. Here we will instead look for solutions that satisfy all of them at once by setting
\begin{equation}
\alpha(K'\mid K_n)\equiv\min\left\{ 1,e^{-\left(\sum_{i=1}^3{\lambda_iz_iD_i}\right)} \right\}
\end{equation}
where $z_i$ are as-of-yet undetermined normalization terms.

\subsection{Motivation}

We first take a moment to examine the sensitivities $\lambda_i$, and impose the following condition: $\lambda_1>\lambda_2$ and $\lambda_1>\lambda_3$. Not doing so would mean the algorithm could give us some false positives. This leads us to notice that a key aspect of the MHMCMC method is information. Due to the ill-posed nature of the problem, we need to keep every piece of information that can be gleaned. We will keep this idea in mind throughout the later sections.\\

As for the normalizations proper, the naive approach to our problem would be to divide each data term by a constant value. In this formulation, our normalization terms would have the form
\begin{equation}
z_i=\frac{1}{c_i} 
\end{equation}
where $c_i$ can be determined by looking at representative values of our data terms.\\

This approach has one advantage, which is that it retains information very well. The relationships between quantities is affected by a constant factor, and its evolution is therefore preserved across iterations. Unfortunately, this method is very unstable, and is not particularly viable. One can think of the opposite method to this one being dividing each data term by itself. Clearly, this would erase all information contained within our results, but it would successfully normalize it, given a broad enough definition of success.\\

Concretely, we seek to find a normalization that delivers information about the evolution of our data terms, but bounds the results so that we may control their magnitudes and work with their relative relationships.

\section{Normalized with Inertia}
\label{sec:norminertia}
We introduce the concept of inertia in this framework. Inertia can be thought of as the weight (call it $w$) being applied to either previous method. Though we do not want to divide by only a constant, there is merit to letting some information trickle through to us. If we do not bound the quantities we are examining, then we will obtain very small or very large values for $\alpha$, effectively $0$ or $1$, which is undoing the work of the MHMCMC. We attempt to bound our likelihood externally. We define $\alpha_h$ such that
\begin{equation}
\alpha(K'\mid K_n)\equiv z_0\alpha_h=z_0e^{-\left(\sum_i{\lambda_iz_iD_i}\right)}
\end{equation}

\subsection{Global Normalizations}

Even a cursory analysis of our early attempts at solving this heat conductivity problem have revealed a desperate need to correctly normalize our data in order to get meaningful likelihoods. Some issues of note have been the idea that the inertia of the process, the value of previous guesses, contains information which is important to the successful convergence of our algorithm. Another is the fact that the variance of data terms means that we require a strong normalization term, at the expense, perhaps, of information, if we are to obtain meaningful results.\\

Addressing the second point, we decide to deviate slightly from one aspect of our method, and use a global result. Computationally, we will only be tracking one variable, and this poses no problem. But note that using a global result in computing $\alpha$ implies that our process is no longer a Markov process, as the probability of reaching the next step is dependent on the past and not just the present.

\subsection{Formulation of $Z^{(1)}$}
\label{sec:z1}
First, let $ \alpha_{h,m}=\max_j\left\{ \alpha_{h,j} \right\}, \ \forall j $ and $D_{i,m}=\max_j\left\{ D_{i,j} \right\}, \ \forall j$. We denote $Z^{(1)}$ the normalization
\begin{align}
	z_{0,j}^{(1)}&=w_0\frac{1}{\alpha_{h,j}}+(1-w_0)\frac{1}{\alpha_{h,m}} \\
	z_{i,j}^{(1)}&=w\frac{1}{\abs{D_{i,j}}}+(1-w)\frac{1}{\abs{D_{i,m}}}
\end{align}
While this effectively bounds our acceptance probability between $[0,1]$, it does so at the expense of the Markov property of our algorithm. Removing this property exhibits some instability in the evolution of the algorithm. Namely, they appear to converge to false positives, an effect which must be explored more fully.

\subsection{Restricted Random Interval}

Examining the values of $\alpha$ that we now produce reveals that we have greatly tightened the spread. Almost all of our values are contained in a narrow band (which changes depending on parameters), say between $0.6$ and $0.75$. Again, this means we are losing information, as the difference in the values of $\alpha$ are lost by comparing them over the entire $[0,1]$ interval.\\

We change the 2nd step in the MHMCMC algorithm, which was $\alpha>u\sim U[0,1]\Rightarrow K_{n+1}=K'$. We now restrict the interval over which we draw $u$, taking its lower and upper bounds at the $j$th iteration to be $[u_{\min},u_{\max}]$, where for some small constant $\zeta$,
\begin{equation}
u_{\min}=\min_{i<j}{\alpha_i}-\zeta \quad\wedge\quad u_{\max}=\max_{i<j}{\alpha_i}+\zeta
\end{equation}
While perhaps more restrictive, this formulation also greatly increases the speed at which the algorithm begins to converge by effectively selecting those guesses which are the most promising, relative to the past performance of the algorithm. This method implies that we will not, with probability $1$, decide the outcome of a guess, they simply become (as per $\zeta$) extremely unlikely to be accepted or rejected.

\section{Locally Focused Normalization}
\label{sec:localnorm}
We now attempt to modify $Z^{(1)}$ in order to retain the original Markov property of the algorithm. The property was violated in the second term, which unfortunately also guarantees we bound our results.

\subsection{Formulation of $Z^{(2)}$}

Denote a new normalization scheme $Z^{(2)}$, given by
\begin{align}
	z_{0,j}^{(2)}&=w_0\frac{1}{\alpha_{h,j}}+(1-w_0)\frac{1}{\alpha_{h,j-1}} \\
	z_{i,j}^{(2)}&=w\frac{1}{\abs{D_{i,j}}}+(1-w)\frac{1}{\abs{D_{i,j-1}}} 
\end{align}
While we have recovered the Markov property, we must now contend with unbounded values for $\alpha$. We note now that preliminary attempts to use $z_{i,j}^{(2)}$ with $z_{0,j}^{(1)}$ did not yield promising results.\\

While this formulation provides good results, it does require us to find an empirical bound for $\alpha$, as it is no longer bounded by $z_0$. For the results presented below, we imposed $\alpha\in[0,1.5]$, setting
\begin{equation}\label{eq:z4alpha}
\alpha\left(K'\mid K_n\right)=\min\left\{ 1.5,z_0e^{-\left(\sum_{i=1}^3{\lambda_iz_iD_i}\right)} \right\}
\end{equation}

\subsection{Results}

\begin{figure}[h]
	\centering
	\subfigure[Reconstruction using $Z^{(1)}$.]{
		\includegraphics[height=3cm]{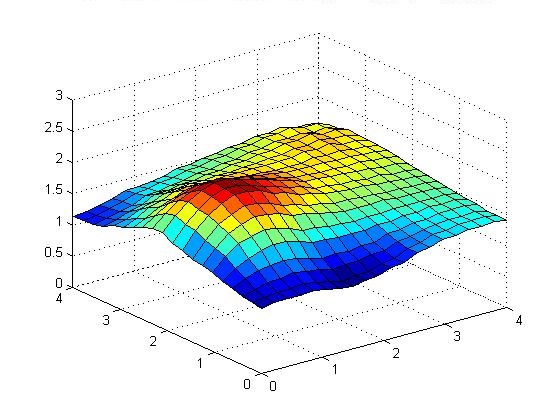}}
	\subfigure[Reconstruction using $Z^{(2)}$.]{
		\includegraphics[height=3cm]{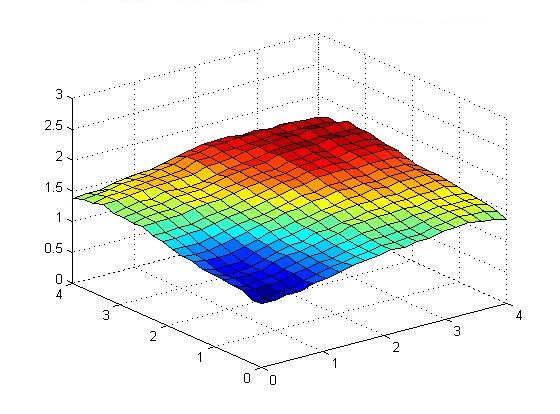}}
	\caption{$Z^{(1)}$ and $Z^{(2)}$ reconstructions of a tilted plane with priors, $2$ million iterations.}
	\label{fig:tiltedz1z4}
\end{figure}

\begin{figure}[h]
	\centering
	\subfigure[Reconstruction using $Z^{(1)}$.]{
		\includegraphics[height=3cm]{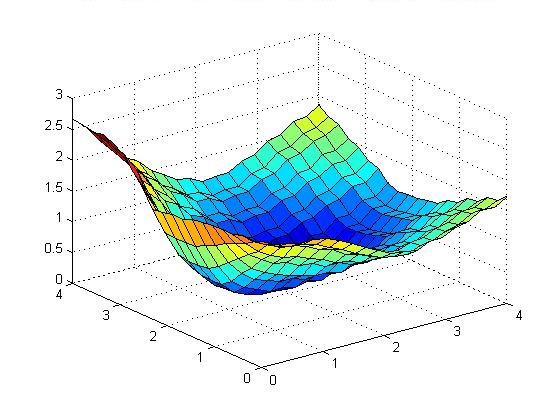}}
	\subfigure[Reconstruction using $Z^{(2)}$.]{
		\includegraphics[height=3cm]{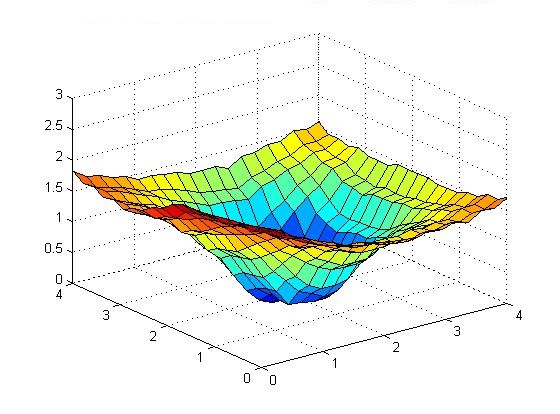}}
	\caption{$Z^{(1)}$ and $Z^{(2)}$ reconstructions of a Gaussian well with priors, $4$ million iterations.}
	\label{fig:gaussz1z4}
\end{figure}

The parameters we have to determine are $\lambda_1,\lambda_2,\lambda_3,w,w_0$ and the cutoff for $\alpha$ as in~\eqref{eq:z4alpha}. We have concluded we must set $\lambda_1>\lambda_i,\ \forall i>1$ and we have by definition $w,w_0\in[0,1]$. The exact values of the sensitivities and inertia factors are at the moment heuristically chosen to be
\begin{align*}
\lambda_1=0.5,\ \lambda_2&=0.15,\ \lambda_3=0.45\\
w_0=0.1,\ w&=0.75,\ \alpha_{\textrm{cutoff}}=1.5
\end{align*}
For the tilted plane, we obtain Figure~\ref{fig:tiltedz1z4}.
As mentioned in Section~\ref{sec:z1}, we have some instability in the form of incorrect convergence for $Z^{(1)}$, which is apparent in Figure~\ref{fig:gaussz1z4} as well. On the other hand, $Z^{(2)}$ converges well and produces a smooth reconstruction. We can also note that it achieves slightly better results than the no-normalizations case in only $2$ million iterations.
The instability in $Z^{(1)}$ is again apparent, and leads us to conclude that the loss of the Markov property in the algorithm may be detrimental to its performance. However, the reconstruction of the Gaussian well has substantially improved when using $Z^{(2)}$. It achieves a smoother reconstruction as without normalizations (see Figure~\ref{fig:gaussfirst}), and in $4M$ iterations instead of $10M$.\\

Going back to our error metric $\Gamma$, we see the improvement manifest itself rather clearly, with acceptances being on the order of $\sim55\%$ instead $\sim95\%$ as they were before.
\begin{figure}[h]
	\centering
	\subfigure[$\Gamma_{Z^{(2)}}$ for tilted plane.]{
		\includegraphics[height=3cm]{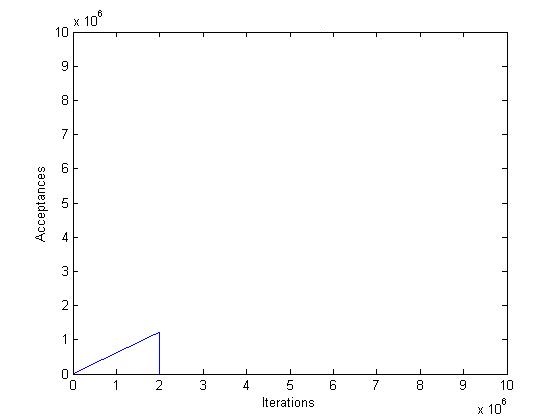}}
	\subfigure[$\Gamma_{Z^{(2)}}$ for Gaussian well.]{
		\includegraphics[height=3cm]{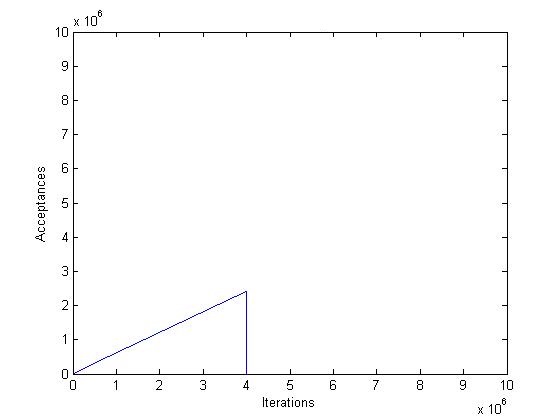}}
	\label{fig:gammaz1z4}
	\caption{Plots of $\Gamma$ for $Z^{(2)}$ reconstructions with priors.}
\end{figure}

\section{Conclusion}

The need for normalizing factors arose from the variance in the magnitudes of data terms $D_i$ from one iteration to the next. In formulating those factors, we focused on conserving the information contained in $D_i$ while bounding our quantities, and we confirmed the importance of retaining the Markov property in this context. However, by using the $Z^{(2)}$ formulation, we were able to obtain faster and better reconstructions of the conductivity for both the tilted plane and the Gaussian well.\\

Despite the encouraging results, several avenues need to be explored more fully. The long-run behavior of $Z^{(2)}$ seems to exhibit some stagnation, seemingly having converged as best as it can. In addition, very preliminary results have been obtained for a scheme that lies between $Z^{(1)}$ and $Z^{(2)}$, which updates the $(1-w)$ terms only when a guess is accepted, has shown competitive performance relative to $Z^{(2)}$.\\

As the algorithm currently stands $\alpha_{\textrm{cutoff}}$, the sensitivities $\lambda_i$, and the inertia factors $w,w_0$ must be determined heuristically. It is possible we may be able to dynamically adjust the values as the algorithm runs, through a constrained optimization of the acceptance rate, but that remains to be studied.\\

Finally, we would like to implement Gaussian-Mixture based MCMC algorithms, that treat the proposal density as an unknown to be approximated, and combine this framework with our normalization schemes to observe the interaction of the two methods.

\ifCLASSOPTIONcaptionsoff
  \newpage
\fi



\bibliographystyle{IEEEtran}

\begin{thebibliography}{9}
	
	\bibitem{fox}
	Fox, C., Nicholls, G., Tan, S.
	\emph{Inverse Problems, Physics 707},
	The University of Auckland, ch. 7-9.
	
	
	\bibitem{hastings}
	Hastings, W.
	``Monte Carlo Sampling Methods Using Markov Chains and Their Applications."
	\emph{Biometrika},
	Vol 57, No. 1, (1970), pp. 97-109.
	
	\bibitem{luengo}
	Luengo, D., Martino, L.
	``Fully Adaptive Gaussian Mixture Metropolis-Hastings Algorithm" 
	\emph{Proc. ICASSP 2013}, Vancouver (Canada), pp. 6148-6152.
	
	\bibitem{metropolis}
	Metropolis, N., Rosenbluth, A., et. al.
	``Equations of State Calculations by Fast Computing Machines" \emph{Journal of Chemical Physics},
	Vol 21 (1953), pp. 1087-1092.
	
	\bibitem{rosenthal}
	Rosenthal, J.
	``Optimal Proposal Distributions and Adaptive MCMC." Chapter for \emph{MCMC Handbook} (2010), avail. at http://www.probability.ca/jeff/ftpdir/galinart.pdf
	
	\bibitem{sauer}
	Sauer, T.,
	\emph{Numerical Analysis},
	Pearson Addison-Wesley, 2006.
	
	\bibitem{zambelli}
	Zambelli, A.,
	``A Multiple Prior Monte Carlo Method for the Backward Heat Diffusion Problem"
	\emph{Proc. CMMSE 2011}, Benidorm (Spain), 
	Vol 3, pp. 1192-1200.
	
\end{thebibliography}
%

%
\newpage

\begin{IEEEbiographynophoto}{Antoine E. Zambelli}
received a Bachelor of Arts in Pure Mathematics from the University of California, Berkeley in 2011 and a Masters in Financial Engineering from the University of California, Los Angeles in 2014. His personal research interests include inverse problems, nonlinear dynamics, global optimizers, and financial derivatives. E-mail: antoine.zambelli@cal.berkeley.edu.
\end{IEEEbiographynophoto}





\end{document}